\documentclass[10pt]{amsart}
\usepackage[T1]{fontenc}
\usepackage{lmodern} 
\usepackage{amssymb}
\usepackage{enumitem}
\usepackage[colorlinks=true,linkcolor=blue,citecolor=blue,urlcolor=blue]{hyperref}

\newtheorem{theorem}{Theorem}[section]
\newtheorem{lemma}[theorem]{Lemma}
\newtheorem{proposition}[theorem]{Proposition}
\newtheorem{corollary}[theorem]{Corollary}
\theoremstyle{definition}
\newtheorem{definition}[theorem]{Definition}
\newtheorem{question}[theorem]{Question}
\theoremstyle{remark}
\newtheorem{remark}[theorem]{Remark}

\newcommand{\Q}{\mathcal{Q}}
\newcommand{\A}{\mathcal{A}}
\newcommand{\B}{\mathcal{B}}
\newcommand{\E}{\mathcal{E}}
\newcommand{\R}{\mathcal{R}}
\newcommand{\N}{\mathbb{N}}
\newcommand{\up}{{\uparrow}}
\newcommand{\dn}{{\downarrow}}
\newcommand{\cl}[1]{\overline{#1}}
\newcommand{\Sig}{\Sigma}
\begin{document}

\title{An algebraic dcpo whose Scott space is not co-sober}

\ifdefined\ANONYMOUS\else
\author{Kia Keng Giam}
\address{Temasek Junior College, Singapore}
\email{giamkiakeng@gmail.com}

\author{Weng Kin Ho}
\address{National Institute of Education, Nanyang Technological University, Singapore}
\email{wengkin.ho@nie.edu.sg}

\author{Bei Liu}
\address{School of Mathematical Sciencs, Hebei Normal University, China}
\email{liubei@hebtu.edu.cn}
\fi


\begin{abstract}
A nonempty compact saturated set is k-irreducible if it is not the union of two proper compact
saturated subsets. A $T_0$ space is co-sober if every k-irreducible set is the saturation of a
point. Zhao, He and Wang proved that every locally compact coherent space is co-sober. We construct an algebraic
dcpo whose Scott space is locally compact and sober, but not co-sober. So their theorem
needs coherence, even for Scott spaces of algebraic domains. Hou and Li showed that quasicontinuous
domains such as algebraic dcpos are determined up to isomorphism by the compact saturated sets of their Scott spaces;
we extend this to posets with locally compact well-filtered Scott spaces. They asked whether every
poset determined in this way has a co-sober Scott space. Our example shows that the answer is no, and that
sobrification does not preserve co-sobriety.
\end{abstract}

\maketitle

\section{Introduction}

A space is sober if every irreducible closed set is the closure of a unique point. Escard\'o,
Lawson and Simpson \cite[Definition~8.9]{ELS04} considered the analogue in which compact saturated
sets take the place of closed sets. A nonempty compact saturated set is \emph{k-irreducible} if it
is not the union of two proper compact saturated subsets. A $T_0$ space is \emph{co-sober} if
every k-irreducible set is the saturation $\up x=\{y:x\le y\}$ of a point $x$. Every set $\up x$ is
k-irreducible, and co-sobriety says that there are no others. For instance, every Hausdorff space
is co-sober: if $p\ne q$ lie in a compact set $Q$ and $U\ni p$, $V\ni q$ are disjoint open sets,
then $Q=(Q\setminus U)\cup(Q\setminus V)$ is a union of two proper compact subsets. The same
authors asked whether every sober space is co-sober
\cite[Problem~9.7]{ELS04}. Wen and Xu \cite{WX18} gave a negative answer, and Shen, Wu, Xi and Zhao
\cite{SWXZ20} gave one with the Scott space of a dcpo. Recently Ji and Xu found further dcpos whose
Scott spaces are not co-sober: a complete Boolean algebra \cite{JX26cba} and a countable dcpo
\cite[Corollary~6.3]{JX26sp}.

On the positive side, Zhao, He and Wang \cite{ZHW23} proved that every locally compact
\emph{coherent} space is co-sober; a space is coherent if the intersection of two compact saturated
sets is compact. Is coherence needed, even for sober spaces, and in particular for the Scott
spaces of continuous or algebraic dcpos? It is.

\begin{theorem}\label{thm:domain}
There is an algebraic dcpo $P$, with a least and a greatest element, whose Scott space $\Sig P$ is
not co-sober. The space $\Sig P$ is compact, locally compact and sober.
\end{theorem}

To our knowledge, Theorem~\ref{thm:domain} gives the first locally compact sober space that is not
co-sober.

Section~\ref{sec:consequences} proves one general result and draws two consequences of
Theorem~\ref{thm:domain}. Let $\Q(X)$ be the set of nonempty compact saturated subsets of a space
$X$, ordered by reverse inclusion. Following Hou and Li \cite[Definition~3.1]{HL25}, a poset $L$ is
\emph{$\Q_\sigma$-unique} if every poset $M$ with $\Q(\Sig M)\cong\Q(\Sig L)$ is isomorphic to $L$;
that is, $L$ is determined by the compact saturated sets of its Scott space. They proved that every
quasicontinuous domain is $\Q_\sigma$-unique \cite[Proposition~3.10]{HL25}; quasicontinuous domains
include all continuous dcpos, for instance all algebraic ones, and have locally compact sober Scott
spaces \cite[Section~III-3]{GHKLMS03}. As their proof uses a result that is cited there as
submitted, we give a complete proof of a more general statement (Proposition~\ref{prop:HL}): every
poset whose Scott space is locally compact and well-filtered is $\Q_\sigma$-unique. Sober spaces are
well-filtered (Section~\ref{sec:prelim}).

Hou and Li asked whether co-sobriety of $\Sig L$ is necessary for $L$ to be $\Q_\sigma$-unique
\cite[Problem~3.22]{HL25}. It is not: $P$ is $\Q_\sigma$-unique, but $\Sig P$ is not co-sober
(Corollary~\ref{cor:HL}). Next, let $Y$ be the poset of compact elements of $P$ with its Alexandroff
topology, whose open sets are the upper sets. Then $\Sig P$ is the sobrification of $Y$, the sober
space with the same lattice of open sets. But $Y$ is co-sober \cite[Lemma~2.1(2)]{JX26sp}. So
sobrification does not preserve co-sobriety (Corollary~\ref{cor:sobrification}).

To prove Theorem~\ref{thm:domain}, we start from the \emph{seed space} $D$ of
Section~\ref{sec:seed}, the space of Xu, Xu, Ji and Qiu \cite[Section~3 and
Proposition~3.3]{XXJQ26k}, which goes back to van Douwen \cite{vD93} and Kov\'ar
\cite[Example~2.2]{Kov05}. Its topology on $\N$ is compact and $T_1$. Certain infinite sets, the
\emph{branches}, each converge to a point, their limit; distinct branches meet in finite sets. The
space $D$ is \emph{anti-Hausdorff}: any two nonempty open sets meet.

The dcpo $P$ is the ideal completion of a base of $D$ ordered by reverse inclusion. Its Scott space
contains a compact saturated set $K$ whose set $\min K$ of minimal elements is a copy of $D$; $K$ is
k-irreducible but, as $\min K$ is not a singleton, not of the form $\up x$. This plan, explained in
\S\ref{ssec:idea}, comes from Xu, Xu, Ji and Qiu \cite[Lemma~4.7 and Theorems~5.1 and~5.3]{XXJQ26s}
and from Ji and Xu \cite[Proposition~3.3 and Theorem~3.4]{JX26cba}. In their examples the copy is
closed in $K$ and is a \emph{KC space} (its compact subsets are closed), so for every compact
$K'\subseteq K$ the trace $K'\cap\min K$ is compact, hence closed in $\min K$. As $\min K$ is
anti-Hausdorff, this makes $K$ k\nobreakdash-irreducible (see the proof of
Theorem~\ref{thm:domain}). In a sober space the copy cannot be closed in $K$
(Remark~\ref{rem:domain}(c)). We obtain closed traces instead from a base smaller than that of all
nonempty open sets, defined by one finiteness rule per branch (Definition~\ref{def:base},
Lemma~\ref{lem:dtrace}).

Moreover, $K$ is not \emph{prime}: it is included in the union of two compact saturated sets but
in neither of them (Remark~\ref{rem:domain}(a)). The dcpo $P$ is not a lattice, and $\Sig P$ is not
first-countable (Remark~\ref{rem:domain}(b) and~(e)). We end with the question whether every
first-countable, locally compact, sober space is co-sober (Question~\ref{q:firstcountable}); the
method of this paper cannot give a first-countable counterexample (Remark~\ref{rem:domain}(e)).

\section{Preliminaries}\label{sec:prelim}

All spaces are $T_0$, and $\N=\{0,1,2,\dots\}$. The \emph{specialization order} of a space $X$ is
$x\le y\iff x\in\cl{\{y\}}$; thus $\cl{\{y\}}=\dn y$, open sets are upper sets, closed sets are lower
sets, and a set is \emph{saturated} (an intersection of open sets) exactly when it is an upper set.
\emph{Compact} means that every open cover has a finite subcover; no separation is assumed. A
nonempty set is compact if and only if every cover of it by a family of open sets that is closed
under finite unions has a member that includes it. The saturation $\up C$ of a compact set $C$ is
compact, because every open set that includes $C$ includes $\up C$. We write $\Q(X)$ for the set of nonempty compact saturated subsets of $X$, and $\min K$ for the
set of minimal elements of $K\in\Q(X)$.

\begin{definition}\label{def:kirr}
A set $K\in\Q(X)$ is \emph{k-irreducible} \cite[Definition~8.9]{ELS04} if $K=K_1\cup K_2$ with $K_1,K_2\in\Q(X)$ implies $K=K_1$ or
$K=K_2$, and \emph{prime} if $K\subseteq K_1\cup K_2$ with $K_1,K_2\in\Q(X)$ implies $K\subseteq K_1$ or
$K\subseteq K_2$. A set of the form $\up a$ is called \emph{principal}. The space $X$ is
\emph{co-sober} if every k-irreducible member of $\Q(X)$ is principal.
\end{definition}

Principal sets are prime, and prime sets are k-irreducible. A nonempty subset $A$ of $X$ is
\emph{irreducible} if $A\subseteq F_1\cup F_2$ with $F_1,F_2$ closed implies $A\subseteq F_1$ or
$A\subseteq F_2$. The space $X$ is \emph{locally compact} if every point has a neighbourhood base of
compact sets. It is \emph{anti-Hausdorff} if any two nonempty open sets intersect; a nonempty space has this property
if and only if it is irreducible as a subset of itself. It is \emph{well-filtered} if for every
nonempty family $\mathcal K\subseteq\Q(X)$ filtered under inclusion and every open set $U$,
$\bigcap\mathcal K\subseteq U$ implies $K\subseteq U$ for some $K\in\mathcal K$. Sober spaces are well-filtered
\cite[Chapter~8]{GL13}. In a well-filtered space the intersection of a nonempty filtered family
$\mathcal K\subseteq\Q(X)$ is compact (see also \cite[Lemma~2.2(iii)]{Giam26}) and nonempty. It is
nonempty, as no member of $\mathcal K$ is included in the open set $\emptyset$. If a family
$\mathcal U$ of open sets covers $\bigcap\mathcal K$, then $\bigcup\mathcal U$ includes some
$K\in\mathcal K$, so finitely many members of $\mathcal U$ cover $K$, hence $\bigcap\mathcal K$.

\begin{lemma}\label{lem:prime}
In a locally compact space every prime compact saturated set is principal.
\end{lemma}

\begin{proof}
Let $K\in\Q(X)$ be prime without a least element. For each $x\in K$ choose $y_x\in K$ with
$x\not\le y_x$, an open set $U_x\ni x$ with $y_x\notin U_x$, and a compact neighbourhood of $x$ inside
$U_x$; its saturation $N_x\subseteq U_x$ is a compact saturated neighbourhood of $x$. Finitely many
interiors cover $K$, so $K\subseteq N_{x_1}\cup\dots\cup N_{x_n}$, and by primeness and induction
$K\subseteq N_{x_i}\subseteq U_{x_i}$ for some $i$, which contradicts $y_{x_i}\in K\setminus U_{x_i}$.
\end{proof}

In a coherent space every k-irreducible set $K$ is prime. Indeed, let $K\subseteq K_1\cup K_2$ with
$K_1,K_2\in\Q(X)$. By coherence, $K\cap K_1$ and $K\cap K_2$ are compact; they are saturated, and
their union is $K$. So $K=K\cap K_i\subseteq K_i$ for some $i$ (if one of the two sets is empty,
$K$ is the other). With
Lemma~\ref{lem:prime} this gives the theorem of Zhao, He and Wang: every locally compact coherent
space is co-sober \cite{ZHW23}.

Directed sets are nonempty, and a \emph{dcpo} is a poset in which every directed subset has a
supremum. A subset $O$ of a poset $P$ is \emph{Scott-open} if it is an upper set and every directed
set that has a supremum in $O$ meets $O$; the Scott-open sets form the \emph{Scott topology}, and
$\Sig P$ denotes $P$ with this topology. For every $x\in P$ the set $\dn x$ is Scott-closed, so
the specialization order of $\Sig P$ is the order of $P$. A monotone map $f$ between posets that
preserves the existing suprema of directed sets is continuous for the Scott topologies: if $E$ is
directed and $f(\sup E)=\sup f(E)$ lies in a Scott-open set $O$, then the directed set $f(E)$
meets $O$. An element $k\in P$ is \emph{compact} if $k\le\sup E$, with $E$ directed, implies $k\le d$ for some
$d\in E$; and $P$ is \emph{algebraic} if every element of $P$ is the supremum of a directed set of
compact elements. An algebraic dcpo is also called an algebraic domain.

\section{The seed space}\label{sec:seed}

This section constructs the space $D$ and proves its properties (Lemma~\ref{lem:seed}). For a family
$\A$ of subsets of $\N$ and $S\subseteq\N$ let $\A(S)=\{a\in\A: a\cap S\text{ is infinite}\}$.

\begin{proposition}\label{prop:mad}
There are an infinite family $\A$ of infinite subsets of $\N$ and a map $\ell\colon\A\to\N$ such
that
\begin{enumerate}[label=\textup{(M\arabic*)},align=right,leftmargin=3.2em]
  \item\label{M1} $a\cap b$ is finite for distinct $a,b\in\A$;
  \item\label{M2} $\A(S)\ne\emptyset$ for every infinite $S\subseteq\N$;
  \item\label{M3} for every $S\subseteq\N$, if $\A(S)$ is infinite, then $\ell[\A(S)]=\N$.
\end{enumerate}
\end{proposition}

A family of infinite subsets of $\N$ with \ref{M1} is called \emph{almost disjoint}.
Condition \ref{M2} says that $\A$ is a maximal such family, as an infinite $S$ can be added to $\A$
exactly when $\A(S)=\emptyset$. In this form the statement is given by Xu, Xu, Ji and
Qiu \cite[Theorem~3.1]{XXJQ26k}, who cite van Douwen's slicing theorem \cite[Theorem~7.1]{vD93}.
Although known, it is proved at the end of this section; it may be skipped.

Let $(\A,\ell)$ be as in Proposition~\ref{prop:mad}. We call the members of $\A$ \emph{branches}
and $\ell(a)$ the \emph{limit} of $a$; a set $B\cup\{\ell(a)\}$ with $a\in\A$ and $B\subseteq a$ is
a \emph{branch piece}. The \emph{seed space} of $(\A,\ell)$ is $D=(\N,\tau)$, in which every branch
converges to its limit:
\[
  \tau=\{U\subseteq\N:\ \ell(a)\in U\text{ implies that }a\setminus U\text{ is finite, for every }
  a\in\A\}.
\]
By \ref{M3}, every closed set $C\ne\N$ of $D$ has $\A(C)$ finite; with \ref{M2}, this makes $D$
compact and anti-Hausdorff (Lemma~\ref{lem:seed}(ii),(iii)).

\begin{lemma}\label{lem:seed}
Let $(\A,\ell)$ satisfy \ref{M1}--\ref{M3}.
\begin{enumerate}[label=(\roman*),align=right]
  \item $\tau$ is a $T_1$ topology. A set $C\subseteq\N$ is closed in $D$ if and only if $\ell(a)\in
        C$ for every $a\in\A(C)$. Branch pieces are closed and compact.
  \item If $C\ne\N$ is closed in $D$, then $\A(C)$ is finite and $C\setminus\bigcup\A(C)$ is finite.
  \item $D$ is compact and anti-Hausdorff.
  \item If $S\subseteq\N$ and $x\in\cl{S}\setminus S$, then there is $a\in\A(S)$ with $\ell(a)=x$.
  \item Every fibre $\ell^{-1}(x)$, $x\in\N$, is infinite.
\end{enumerate}
\end{lemma}

\begin{proof}
(i) $\emptyset,\N\in\tau$. If $U_i\in\tau$ and $\ell(a)\in\bigcup_iU_i$, then $\ell(a)\in U_j$ for
some $j$, and $a\setminus\bigcup_iU_i\subseteq a\setminus U_j$ is finite. If $U,V\in\tau$ and
$\ell(a)\in U\cap V$, then $a\setminus(U\cap V)=(a\setminus U)\cup(a\setminus V)$ is finite. So $\tau$
is a topology. A set $C$ is closed if and only if, for every $a$, $\ell(a)\notin C$ implies that
$a\cap C$ is finite; this is the stated criterion. Finite sets satisfy it, so $D$ is $T_1$. If
$b\in\A$ meets $B\cup\{\ell(a)\}$ in an infinite set, then $b\cap a$ is infinite, so $b=a$ by
\ref{M1}, and $\ell(b)=\ell(a)$ lies in the branch piece; so branch pieces are closed. Given an open
cover of $B\cup\{\ell(a)\}$, one member $U$ contains $\ell(a)$; then $B\setminus U\subseteq a\setminus
U$ is finite, and finitely many further members cover it.

(ii) Let $C$ be closed. If $\A(C)$ is infinite, then for every $n\in\N$ there is $a\in\A(C)$ with
$\ell(a)=n$ by \ref{M3}, and $n=\ell(a)\in C$ by (i); so $C=\N$. Let $\A(C)$ be finite and
$R=C\setminus\bigcup\A(C)$. If $R$ were infinite, \ref{M2} would give $b\in\A$ with $b\cap R$ infinite;
then $b\in\A(C)$ and $b\cap R=\emptyset$, a contradiction.

(iii) Let $C\ne\N$ be closed. By (ii) and (i), $C$ is the union of the finite set
$C\setminus\bigcup\A(C)$ and of the finitely many branch pieces $(a\cap C)\cup\{\ell(a)\}$,
$a\in\A(C)$; so $C$ is compact. Given an open cover of $\N$, choose a nonempty member $U$; the closed
set $\N\setminus U$ is compact, hence covered by finitely many members. So $D$ is compact. If $U,V$
were disjoint nonempty open sets, then $\N=(\N\setminus U)\cup(\N\setminus V)$ would be a union of two
closed sets different from $\N$, so by (ii) there would be a finite $\mathcal F\subseteq\A$ and a
finite $F\subseteq\N$ with $\N=F\cup\bigcup\mathcal F$. For $a\in\A\setminus\mathcal F$ the set
$a\cap(F\cup\bigcup\mathcal F)$ is finite by \ref{M1}, which is absurd because $a$ is infinite.

(iv) Suppose that no $a\in\A(S)$ has $\ell(a)=x$. Then $\A(S)$ is finite by \ref{M3}, and
$R=S\setminus\bigcup\A(S)$ is finite by the argument of (ii). The set
$R\cup\bigcup_{a\in\A(S)}\bigl((a\cap S)\cup\{\ell(a)\}\bigr)$ is closed by (i), it includes $S$, and
it does not contain $x$, because $x\notin S$ and $x\neq\ell(a)$ for $a\in\A(S)$. This contradicts
$x\in\cl{S}$.

(v) Suppose $\ell^{-1}(x)=\{b_1,\dots,b_r\}$ with $r\ge0$, and let $S=\N\setminus(b_1\cup\dots\cup
b_r)$. Every $a\in\A\setminus\{b_1,\dots,b_r\}$ meets $b_1\cup\dots\cup b_r$ in a finite set, so
$a\in\A(S)$. Hence $\A(S)$ is infinite, and by \ref{M3} some $a\in\A(S)$ has $\ell(a)=x$; then $a=b_i$
for some $i$, but $b_i\cap S=\emptyset$.
\end{proof}

\begin{proof}[Proof of Proposition~\ref{prop:mad}]
Step~1 below is a standard lemma on almost disjoint families, and Step~2 is an argument of Simon
\cite{Sim80}, which we iterate in Step~3; see Hru\v{s}\'ak's survey \cite[Proposition~2 and
Theorem~46]{Hru14} for both.
Let $\E$ be an infinite maximal almost disjoint family of infinite subsets of $\N$; it exists by
Zorn's lemma, applied to the almost disjoint families that include a fixed partition of $\N$
into infinitely many infinite sets. By maximality, $\E(S)\neq\emptyset$ for every infinite
$S\subseteq\N$. Call $S\subseteq\N$ \emph{positive} if $\E(S)$ is infinite. For a positive set
$H$, call a family $\R\subseteq\E$ \emph{dense below} $H$ if $\R\cap\E(J)\ne\emptyset$ for every
positive $J\subseteq H$. Step~0 collects elementary facts on positive sets and on density.
Steps~1--3 produce a positive set $H$ and pairwise disjoint families
$\E_0,\E_1,\dots\subseteq\E$ that are all dense below $H$. Step~4 then identifies $H$ with $\N$ and
takes the infinite sets $e\cap H$, $e\in\E$, as the members of $\A$; the map $\ell$ sends those
with $e\in\E_n$ to $n$.

\emph{Step 0.} If $S$ is positive and $e_1,\dots,e_k\in\E$, then $S'=S\setminus(e_1\cup\dots\cup
e_k)$ satisfies $\E(S')=\E(S)\setminus\{e_1,\dots,e_k\}$, because every other member of $\E$ meets
$e_1\cup\dots\cup e_k$ in a finite set; in particular $S'$ is positive. Consequently, if $\R$ is
dense below $H$, then $\R\cap\E(J)$ is infinite for every positive $J\subseteq H$: if it were a
finite set $\{e_1,\dots,e_k\}$, the positive set $J\setminus(e_1\cup\dots\cup e_k)$ would contradict
density. Moreover, if $\R$ is dense below $H$ and $H'$ is a positive set with $H'\setminus H$ finite,
then $\R$ is dense below $H'$: for a positive $J\subseteq H'$ the set $J\cap H$ differs from $J$ by
finitely many points, so it is positive and $\E(J\cap H)=\E(J)$.

\emph{Step 1. If $H_0\supseteq H_1\supseteq\cdots$ are positive, there is a positive $H\subseteq H_0$
with $H\setminus H_n$ finite for all $n$.} First, $\mathcal F=\bigcap_n\E(H_n)$ is infinite.
Otherwise let $H_n'=H_n\setminus\bigcup\mathcal F$; by Step~0 these sets are positive, they decrease,
and $\bigcap_n\E(H_n')=\emptyset$. Choose distinct points $p_n\in H_n'$. Some $e\in\E$ meets
$\{p_n:n\in\N\}$ in an infinite set; as $p_k\in H_n'$ for $k\ge n$, this $e$ belongs to every
$\E(H_n')$, a contradiction. Now choose distinct $e_0,e_1,\dots\in\mathcal F$. For each $n$ the sets
$e_n\cap H_k$, $k\in\N$, are infinite and decrease, so there are points $y_{n,0}<y_{n,1}<\cdots$
with $y_{n,j}\in e_n\cap H_{\max(n,j)}$. Put $J_n=\{y_{n,j}:j\in\N\}$ and $H=\bigcup_nJ_n$. Then
$J_n\subseteq H_n\subseteq H_0$, and $H\setminus H_k\subseteq\{y_{n,j}:n<k,\ j<k\}$ is finite. As
$e_n\cap H\supseteq J_n$ is infinite for every $n$, the set $H$ is positive.

\emph{Step 2 \textup{(Simon)}. If $\R$ is dense below $H$, there are a positive $H'\subseteq H$ and a
partition of $\R$ into two sets that are both dense below $H'$.} Suppose not: for every positive
$H'\subseteq H$ and every partition of $\R$ into two sets, one of the two sets is not dense below
$H'$. For $e\in\R$ let $f_e\in\{0,1\}^{\N}$ be the characteristic function of $e$, and let
$\R_n^i=\{e\in\R:f_e(n)=i\}$ for $n\in\N$,
$i\in\{0,1\}$. Put $H_{-1}=H$. Given a positive $H_{n-1}\subseteq H$, the sets $\R_n^0$ and $\R_n^1$
partition $\R$, so one of them, say $\R_n^{g(n)}$, is not dense below $H_{n-1}$: there is a positive
$H_n\subseteq H_{n-1}$ with $\R_n^{g(n)}\cap\E(H_n)=\emptyset$. By Step~1 there is a positive
$H_\infty\subseteq H$ with $H_\infty\setminus H_n$ finite for all $n$. At most one member of $\R$
lies in $\E(H_\infty)$: if $e\in\R\cap\E(H_\infty)$, then $e\in\E(H_n)$, so $f_e(n)\neq g(n)$, for
all $n$; that is, $f_e=1-g$, which determines $e$. But $\R\cap\E(H_\infty)$ is infinite by Step~0.

\emph{Step 3. There are a positive set $H$ and pairwise disjoint families $\E_0,\E_1,\dots\subseteq\E$
that are all dense below $H$.} The family $\R_0=\E$ is dense below the positive set $H_0=\N$.
Applying Step~2 repeatedly we get positive sets $H_0\supseteq H_1\supseteq\cdots$ and partitions
$\R_n=\E_n\cup\R_{n+1}$ into two families that are dense below $H_{n+1}$. The families $\E_n$ are
pairwise disjoint, and $\E_n$ is dense below $H_m$ for all $m>n$. Take $H$ from Step~1; by Step~0
every $\E_n$ is dense below $H$.

\emph{Step 4.} Fix a bijection $\theta\colon\N\to H$. Let $\A=\{e\cap H:e\in\E(H)\}$, and put
$\ell(e\cap H)=\theta(n)$ if $e\in\E_n$ and $\ell(e\cap H)=\theta(0)$ if $e$ belongs to no $\E_n$;
this is well defined, because distinct members of $\E(H)$ meet $H$ in distinct sets. Then $\A$ is
an infinite almost disjoint family of infinite subsets of $H$. For an infinite $S\subseteq H$, a set
$e\in\E$ belongs to $\E(S)$ if and only if $e\in\E(H)$ and $e\cap H\in\A(S)$. So $\A(S)\ne\emptyset$;
and if $\A(S)$ is infinite, then $S$ is positive and, for every $n$, some $e\in\E_n$ belongs to
$\E(S)$, that is, $\theta(n)\in\ell[\A(S)]$. Thus $(\A,\ell)$ satisfies \ref{M1}--\ref{M3} on the
countably infinite set $H$, and we transport it to $\N$ along $\theta^{-1}$.
\end{proof}

\section{The construction}\label{sec:construction}

\subsection{The idea}\label{ssec:idea}
\emph{The plan.} We order a base $\B$ of $D$ by reverse inclusion; the points of its ideal
completion $P$, called ideals, are the upward closed and downward directed subsets of
$(\B,\subseteq)$. The basic neighbourhoods of $m\in\N$ form a point $I_m$ of $P$, the map
$m\mapsto I_m$ embeds $D$ in $\Sig P$, and the saturation $K$ of the image is compact and has no
least element (Lemma~\ref{lem:embed}). The \emph{trace} of a compact set $K'\subseteq K$ is
$\{m\in\N:I_m\in K'\}$. If all traces are closed, then $K$ is k-irreducible, because $D$ is
anti-Hausdorff (see the proof of Theorem~\ref{thm:domain}).

\emph{The obstacle.} We must show that every trace is closed (Lemma~\ref{lem:dtrace}), that is, by
Lemma~\ref{lem:seed}(i): if $K'$ contains $I_d$ for infinitely many points $d$ of a branch $a$ with
limit $x$, then $I_x\in K'$. Suppose not. For a contradiction, it suffices to give each point $y$ of
$K'$ a neighbourhood that contains only finitely many of these $I_d$, because finitely many such
neighbourhoods cover the compact set $K'$. This is easy if $y$ lies above some $I_m$ with $m\ne x$
(Cases~1 and~2); so the main step concerns the points $y$ strictly above $I_x$. If
$W\in y$, the set $[W]$ of ideals that contain $W$ is a neighbourhood of $y$ that contains $I_d$
only if $d\in W$. So it suffices to find $W\in y$ with $W\cap a$ finite. For the base of all
nonempty open sets the plan fails (Remark~\ref{rem:domain}(d)).

\emph{The remedy.} We pair the branches by a map $c\colon\A\to\A$ with $c(c(a))=a$ and
$\ell(c(a))\ne\ell(a)$ (Lemma~\ref{lem:pairing}); the branch piece $Z_a$ of $c(a)$ defined in
\S\ref{ssec:base} does not contain $\ell(a)$. The thinned base $\B$ consists of the nonempty open
sets $W$ that obey, for each branch $a$, the rule $(\mathrm R_a)$ of Definition~\ref{def:base}: if
$W$ misses $\ell(a)$ and $Z_a$, then $W\cap a$ is finite. The condition on $Z_a$ is needed. Without
it, each $W\in\B$ would meet every branch $a$ with $\ell(a)\notin W$ in a finite set; so $W$ would
be closed (Lemma~\ref{lem:seed}(i)), hence equal to $\N$, as $D$ is anti-Hausdorff. With the rules,
each ideal $y$ strictly above $I_x$, where $x=\ell(a)$, contains some $W\in\B$ that misses $x$ and
$Z_a$ (because $\N\setminus Z_a$ is a basic neighbourhood of $x$), hence with $W\cap a$ finite
(Lemma~\ref{lem:dtrace}, Case~3).

\emph{The price} is the proof that $\B$ is a base (Lemma~\ref{lem:base}): to turn an open set
$G\ni m$ into a member of $\B$, we cut out of $G$ the branches whose rules apply, and the pairing keeps
these finitely many, as they lie in $\A(\N\setminus G)\cup c[\A(\N\setminus G)]$.

\subsection{Pairing the branches}
The pairing $c$ needed in \S\ref{ssec:idea} maps each fibre $\ell^{-1}(x)$ injectively into its
complement, so the inequality
\eqref{eq:normal} below is necessary; Lemma~\ref{lem:normal} arranges it, and
Lemma~\ref{lem:pairing} shows that it suffices.

\begin{lemma}\label{lem:normal}
There is a pair $(\A,\ell)$ satisfying \ref{M1}--\ref{M3} such that
\begin{equation}\label{eq:normal}
  |\ell^{-1}(x)|\le|\A\setminus\ell^{-1}(x)|\qquad\text{for every }x\in\N .
\end{equation}
\end{lemma}

\begin{proof}
Let $(\A,\ell)$ satisfy \ref{M1}--\ref{M3}. If \eqref{eq:normal} holds, there is nothing to prove;
suppose that it fails for a point $x_0$. Then $\mathcal F=\ell^{-1}(x_0)$ satisfies $|\mathcal F|=|\A|$, because the infinite set $\A$
is the union of $\mathcal F$ and the smaller set $\A\setminus\mathcal F$. We delete $x_0$ and split
$\mathcal F$ between two other points: put $N'=\N\setminus\{x_0\}$
and $a'=a\setminus\{x_0\}$ for $a\in\A$; the map $a\mapsto a'$ is injective by \ref{M1}. Fix distinct
points $x_1,x_2\in N'$ and a partition of $\mathcal F$ into two sets $\mathcal F_1,\mathcal F_2$ of
cardinality $|\mathcal F|$. Let $\A'=\{a':a\in\A\}$ and
\[
  \ell'(a')=\begin{cases}\ell(a)&\text{if }\ell(a)\ne x_0,\\ x_i&\text{if }a\in\mathcal F_i,\ i=1,2.
  \end{cases}
\]
The pair $(\A',\ell')$ satisfies \ref{M1} and \ref{M2} on the countably infinite set $N'$. For
\ref{M3} let $S\subseteq N'$ with $\A'(S)$ infinite, and let $n\in N'$. Then $\A(S)$ is infinite, so
some $a\in\A(S)$ has $\ell(a)=n$; as $n\ne x_0$, we get $\ell'(a')=n$, and $a'\in\A'(S)$. Finally,
for every $x\in N'$ the set $\A'\setminus\ell'^{-1}(x)$ includes $\{a':a\in\mathcal F_1\}$ or
$\{a':a\in\mathcal F_2\}$, so it has cardinality $|\A'|$, and \eqref{eq:normal} holds for
$(\A',\ell')$. Transport this pair to $\N$ along a bijection.
\end{proof}

\begin{lemma}\label{lem:pairing}
Let $(\A,\ell)$ be as in Lemma~\ref{lem:normal}. Then there is a map $c\colon\A\to\A$ such that
$c(c(a))=a$ and $\ell(c(a))\ne\ell(a)$ for every $a\in\A$.
\end{lemma}

\begin{proof}
Write $\kappa=|\A|$ and $\kappa_x=|\ell^{-1}(x)|$, an infinite cardinal by Lemma~\ref{lem:seed}(v).
Then $\kappa=\sup_x\kappa_x$, because $\A$ is the union of countably many fibres, each of them
infinite. We first find a set $E\subseteq\N$ with
$|\ell^{-1}(E)|=|\ell^{-1}(\N\setminus E)|=\kappa$.
If $\kappa_{x_0}=\kappa$ for some $x_0$, take $E=\{x_0\}$; then
$|\ell^{-1}(\N\setminus E)|=|\A\setminus\ell^{-1}(x_0)|\ge\kappa_{x_0}=\kappa$ by \eqref{eq:normal}.
Otherwise $\kappa_x<\kappa$ for all $x$. Enumerate $\N=\{z_1,z_2,\dots\}$ and choose points
$x_1,x_2,\dots$ one after the other with $\kappa_{x_{i+1}}>\max(\kappa_{x_i},\kappa_{z_i})$; this is
possible because no $\kappa_x$ equals the supremum $\kappa$. The points $x_i$ are distinct, and
$\sup_i\kappa_{x_{2i}}=\sup_i\kappa_{x_{2i+1}}=\kappa$. Take $E=\{x_{2i}:i\ge1\}$; then
$|\ell^{-1}(E)|\ge\sup_i\kappa_{x_{2i}}=\kappa$ and
$|\ell^{-1}(\N\setminus E)|\ge\sup_i\kappa_{x_{2i+1}}=\kappa$.

Now let $\beta\colon\ell^{-1}(E)\to\ell^{-1}(\N\setminus E)$ be a bijection and put
$c=\beta\cup\beta^{-1}$. Then $c\circ c$ is the identity of $\A$, and exactly one of $\ell(a)$,
$\ell(c(a))$ lies in $E$.
\end{proof}

\subsection{The thinned base}\label{ssec:base}
From now on $(\A,\ell)$ and $c$ are as in Lemma~\ref{lem:pairing}, and $D=(\N,\tau)$ is the seed
space of $(\A,\ell)$; Lemma~\ref{lem:seed} applies. For $a\in\A$ let
$Z_a=\bigl(c(a)\setminus\{\ell(a)\}\bigr)\cup\{\ell(c(a))\}$ and $O_a=\N\setminus Z_a$. As a
branch piece of $c(a)$, $Z_a$ is closed, and $\ell(a)\notin Z_a$ because $\ell(c(a))\ne\ell(a)$; so
$O_a$ is an open neighbourhood of $\ell(a)$.

\begin{definition}\label{def:base}
$\B$ is the set of all nonempty open subsets $W$ of $D$ that satisfy, for every $a\in\A$, the rule
\[
  (\mathrm R_a)\qquad\text{if }\ell(a)\notin W\text{ and }W\cap Z_a=\emptyset,\text{ then }W\cap a
  \text{ is finite.}
\]
\end{definition}

We say that $(\mathrm R_a)$ \emph{applies} to $W$ if its hypothesis holds: $\ell(a)\notin W$ and
$W\subseteq O_a$.

\begin{lemma}\label{lem:base}
The family $\B$ has the following properties.
\begin{enumerate}[label=(\roman*),align=right]
  \item $\N\in\B$; $O_a\in\B$ for every $a\in\A$; and $\N\setminus\{x\}\in\B$ for every $x\in\N$.
  \item $\B$ is a base of $D$.
\end{enumerate}
\end{lemma}

\begin{proof}
(i) We show that no rule applies to these nonempty open sets. No rule applies to $\N$, and none to
$\N\setminus\{x\}$, because no infinite set $Z_a$ is disjoint from $\N\setminus\{x\}$. Let
$a,b\in\A$ with $\ell(b)\notin O_a$ and $O_a\cap Z_b=\emptyset$. Then $Z_b\subseteq Z_a$, so the
infinite set $c(b)\setminus\{\ell(b)\}$ is included in $c(a)\cup\{\ell(c(a))\}$; hence $c(b)\cap
c(a)$ is infinite, $c(b)=c(a)$ by \ref{M1}, and $b=a$ because $c$ is injective. But $\ell(a)\in
O_a$. So no rule applies to $O_a$.

(ii) Let $m\in G$ with $G$ open; we look for $W\in\B$ with $m\in W\subseteq G$. We cut out of $G$,
except for the point $m$, the branches whose rules apply, and repeat this until no new rule applies.
Put $W_0=G$ and define recursively
\begin{gather*}
  C_j=\N\setminus W_j,\qquad
  T_j=\{a\in\A:\ \ell(a)\notin W_j\text{ and }W_j\cap Z_a=\emptyset\},\\
  W_{j+1}=W_j\setminus\bigcup_{a\in T_j}(a\setminus\{m\});
\end{gather*}
we show along the way that each $W_j$ is open and contains $m$.

\emph{$T_j$ is finite.} The set $C_j$ is closed and $m\notin C_j$. If $a\in T_j$, then $Z_a\subseteq
C_j$, so $c(a)\in\A(C_j)$, which is a finite set by Lemma~\ref{lem:seed}(ii); and $c$ is injective.

\emph{$W_{j+1}$ is open and contains $m$.} Its complement is $C_{j+1}=C_j\cup\bigcup_{a\in
T_j}(a\setminus\{m\})$. Let $b\in\A(C_{j+1})$. As $T_j$ is finite, either $b\cap C_j$ is infinite, and
then $\ell(b)\in C_j$ because $C_j$ is closed; or $b\cap a$ is infinite for some $a\in T_j$, and then
$b=a$ and $\ell(b)=\ell(a)\in C_j$ by the definition of $T_j$. By Lemma~\ref{lem:seed}(i), $C_{j+1}$
is closed. Clearly $m\in W_{j+1}$.

\emph{The sets $T_j$ increase and stay inside one finite set.} $T_j\subseteq T_{j+1}$ because
$W_{j+1}\subseteq W_j$, and $\A(C_{j+1})\subseteq\A(C_j)\cup T_j$ by the previous step. Let
$\mathcal S=\A(C_0)\cup c[\A(C_0)]$, a finite set with $c[\mathcal S]=\mathcal S$, because $c\circ c$ is the identity. By induction on $j$
we get $\A(C_j)\subseteq\mathcal S$ and $T_j\subseteq\mathcal S$: if $a\in T_j$, then
$c(a)\in\A(C_j)\subseteq\mathcal S$, so $a=c(c(a))\in\mathcal S$; and then
$\A(C_{j+1})\subseteq\A(C_j)\cup T_j\subseteq\mathcal S$.

So there is $j$ with $T_{j+1}=T_j$. Put $W=W_{j+1}$. Then $W$ is open and $m\in W\subseteq G$. If
$a\in\A$ satisfies $\ell(a)\notin W$ and $W\cap Z_a=\emptyset$, then $a\in T_{j+1}=T_j$, so $W\cap
a\subseteq\{m\}$. Hence $W\in\B$.
\end{proof}

\subsection{The dcpo and the compact set}

Order $\B$ by reverse inclusion. An \emph{ideal} of $\B$ is a nonempty set $y\subseteq\B$ such that
$U\in y$ whenever $V\in y$, $U\in\B$ and $V\subseteq U$, and for all $U,V\in y$ there is $W\in y$
with $W\subseteq U\cap V$. Let $P$ be the set of all ideals of $\B$, ordered by inclusion, that is,
the ideal completion of $(\B,\supseteq)$. As $\B$ is ordered by reverse inclusion, the principal
ideal of $U\in\B$ consists of the supersets of $U$: it is $\dn U=\{V\in\B:U\subseteq V\}\in P$. Let $[U]=\{y\in P:U\in y\}=\up_P(\dn U)$.

\begin{lemma}\label{lem:idl}
The poset $P$ has the following properties.
\begin{enumerate}[label=(\roman*),align=right]
  \item $P$ is a dcpo in which directed suprema are unions. It is algebraic: its compact elements are
        the principal ideals $\dn U$, $U\in\B$, and every $y\in P$ is the directed union of the
        ideals $\dn U$, $U\in y$.
  \item The sets $[U]$, $U\in\B$, form a base of the Scott topology of $P$, and the specialization
        order of $\Sig P$ is inclusion.
  \item $\Sig P$ is locally compact and sober.
  \item $P$ has the least element $\{\N\}$ and the greatest element $\B$; so $\Sig P$ is compact.
\end{enumerate}
\end{lemma}

\begin{proof}
These are standard facts on ideal completions \cite{GL13}; we include short proofs.

(i) Directed unions of ideals are ideals. If such a union includes $\dn U$, then one of its members
contains $U$ and so includes $\dn U$; thus $\dn U$ is compact. Each $y\in P$ is the directed union of
the ideals $\dn U$, $U\in y$, so a compact $y$ equals one of them.

(ii) By (i), $[U]$ is Scott-open, and a Scott-open set containing $y$ contains some $\dn U$ with
$U\in y$, hence $[U]$; the specialization order is inclusion by Section~\ref{sec:prelim}.

(iii) Each $[U]$ is compact, as it has the least element $\dn U$; with (ii), $\Sig P$ is locally
compact. Let $C$ be closed and irreducible, and $y_C=\{U\in\B:\dn U\in C\}$; it is nonempty and
closed under supersets, as $C$ is a nonempty lower set. For $U,V\in y_C$, some $y\in C$ lies in
$[U]\cap[V]$, as $C$ is irreducible; it contains some $W\subseteq U\cap V$, and $W\in y_C$ as $\dn
W\subseteq y$. So $y_C\in P$, and $y_C\in C$, as it is the directed union of the ideals $\dn U\in C$
and $C$ is Scott-closed. As $y\subseteq y_C$ for $y\in C$, $C$ is the closure of $y_C$.

(iv) Every ideal contains $\N$, the largest member of $\B$. The set $\B$ is an ideal, because any two
nonempty open subsets of $D$ intersect and $\B$ is a base. A space whose specialization order has a
least element is compact.
\end{proof}

For $m\in\N$ let $I_m=\{U\in\B:m\in U\}$.

\begin{lemma}\label{lem:embed}
The sets $I_m$ have the following properties.
\begin{enumerate}[label=(\roman*),align=right]
  \item $I_m\in P$ for every $m$, and the map $e\colon D\to\Sig P$, $e(m)=I_m$, satisfies
        $e^{-1}([U])=U$ for all $U\in\B$.
  \item $I_m\subseteq I_{m'}$ only if $m=m'$. Hence $e$ is a homeomorphism of $D$ onto the subspace
        $e(\N)$ of $\Sig P$.
  \item $K=\up_Pe(\N)=\{y\in P:I_m\subseteq y\text{ for some }m\in\N\}$ is a compact saturated subset
        of $\Sig P$, and $\min K=e(\N)$. In particular $K$ has no least element.
\end{enumerate}
\end{lemma}

\begin{proof}
(i) $\N\in I_m$; $I_m$ is closed under supersets in $\B$; and if $U,V\in I_m$, then $m\in U\cap V$,
which is open, so there is $W\in\B$ with $m\in W\subseteq U\cap V$ by Lemma~\ref{lem:base}(ii), and
$W\in I_m$. Moreover $e(m)\in[U]$ if and only if $m\in U$. (ii) If $m\neq m'$, then
$\N\setminus\{m'\}\in I_m\setminus I_{m'}$ by Lemma~\ref{lem:base}(i). So $e$ is injective; by (i) it
is continuous and $e(U)=[U]\cap e(\N)$ for $U\in\B$, and $\B$ is a base of $D$; so $e$ is a
homeomorphism onto its image. (iii) $e(\N)$ is compact by (ii)
and Lemma~\ref{lem:seed}(iii), so its saturation $K$ is compact. Every element of $K$ is above some
$I_m$, and by (ii) the elements $I_m$ are pairwise incomparable; so they are the minimal elements of
$K$.
\end{proof}

\begin{lemma}\label{lem:dtrace}
Let $K'\subseteq K$ be a compact subset of $\Sig P$. Then its trace $D'=\{m\in\N:I_m\in K'\}$ is
closed in $D$.
\end{lemma}

\begin{proof}
Suppose that there is $x\in\cl{D'}\setminus D'$. By Lemma~\ref{lem:seed}(iv) there is $a\in\A$ with
$\ell(a)=x$ and $a\cap D'$ infinite. Put $S=a\cap D'$ and $C=S\cup\{x\}$; note that $x\notin S$. The
set $C$ and the sets $C\setminus\{d\}$, $d\in S$, are branch pieces, hence closed. As the compact
set $K'$ contains the infinitely many $I_d$, $d\in S$, it suffices to show that every $y\in K'$ has
an open neighbourhood in $\Sig P$ that contains $I_d$ for at most finitely many $d\in S$.

Let $y\in K'$. As $y\in K$, there is $m\in\N$ with $I_m\subseteq y$.

\emph{Case 1: $m\notin C$.} Choose $U\in\B$ with $m\in U\subseteq\N\setminus C$. Then $U\in
I_m\subseteq y$, so $y\in[U]$; and $I_d\notin[U]$ for every $d\in S$, because $d\notin U$.

\emph{Case 2: $m=d\in S$.} Choose $U\in\B$ with $d\in U\subseteq\N\setminus(C\setminus\{d\})$. Then
$y\in[U]$, and for $d'\in S$ we have $I_{d'}\in[U]$ only if $d'=d$.

\emph{Case 3: $m=x$.} As $x\notin D'$, we have $I_x\notin K'$, so $y\neq I_x$ and there is $U\in y$
with $x\notin U$. Also $O_a\in I_x\subseteq y$, by Lemma~\ref{lem:base}(i) and $x=\ell(a)\in O_a$. As
$y$ is an ideal, there is $W\in y$ with $W\subseteq U\cap O_a$. Then $\ell(a)=x\notin W$ and $W\cap
Z_a=\emptyset$, so $W\cap a$ is finite by the rule $(\mathrm R_a)$. Hence $y\in[W]$, and $I_d\in[W]$
only for the finitely many $d\in S\cap W$.
\end{proof}

\begin{proof}[Proof of Theorem~\ref{thm:domain}]
Let $P$ and $K$ be as above. By Lemma~\ref{lem:idl}, $P$ is an algebraic dcpo with a least and a
greatest element, and $\Sig P$ is compact, locally compact and sober. By Lemma~\ref{lem:embed},
$K\in\Q(\Sig P)$ and $K$ is not principal. Let $K=K_1\cup K_2$ with $K_1,K_2\in\Q(\Sig P)$, and let
$D_i=\{m:I_m\in K_i\}$ be the trace of $K_i$. By Lemma~\ref{lem:dtrace}, the traces $D_1$ and $D_2$
are closed in $D$, and $D_1\cup D_2=\N$ because every $I_m$ lies in $K$. As $D$ is anti-Hausdorff, $D_1=\N$ or $D_2=\N$, say
$D_1=\N$. Then $e(\N)\subseteq K_1$, and $K=\up_Pe(\N)\subseteq K_1$ because $K_1$ is saturated. So
$K$ is k-irreducible, and $\Sig P$ is not co-sober.
\end{proof}

\section{Consequences and remarks}\label{sec:consequences}

\subsection{A question of Hou and Li}
In $\Q(X)$ we write $Q\sqsubseteq Q'$ for $Q\supseteq Q'$. Hou and Li also admit the greatest
element $\emptyset$; as order isomorphisms fix it, this does not change $\Q_\sigma$-uniqueness.

\begin{proposition}\label{prop:HL}
Every poset whose Scott space is locally compact and well-filtered is $\Q_\sigma$-unique.
\end{proposition}

For quasicontinuous domains, this is due to Hou and Li \cite[Proposition~3.10]{HL25}.

To prove Proposition~\ref{prop:HL} for a poset $L$, we show that an order isomorphism $F$ of $(\Q(\Sig L),\sqsubseteq)$ onto $(\Q(\Sig M),\sqsubseteq)$, for a poset $M$, maps
principal sets onto principal sets, using two notions defined by the order alone. An element $p$ of
a poset $(Z,\sqsubseteq)$ is
\begin{itemize}[nosep]
  \item \emph{prime} if $u\wedge v\sqsubseteq p$ implies $u\sqsubseteq p$ or $v\sqsubseteq p$,
        whenever the meet $u\wedge v$ exists;
  \item \emph{strongly prime} if $p\in\mathcal K$ whenever $\mathcal K\subseteq Z$ is nonempty,
        compact and saturated in the Scott topology of $(Z,\sqsubseteq)$, and has a meet
        $\bigwedge\mathcal K\sqsubseteq p$.
\end{itemize}
Strong primeness is primeness for Scott-compact families instead of pairs: as $\mathcal K$ is an
upper set, $p\in\mathcal K$ means $k\sqsubseteq p$ for some $k\in\mathcal K$. The second notion is
adapted from \cite[Definition~3.6]{HL25}. In $(\Q(X),\sqsubseteq)$ the meet of two members is their
union, so its prime elements are its prime members in the sense of Definition~\ref{def:kirr}.

\begin{lemma}\label{lem:sprime}
Strongly prime elements and principal sets are related as follows.
\begin{enumerate}[label=(\roman*),align=right]
  \item \textup{(essentially Hou and Li \cite[Lemma~3.18]{HL25})} For every poset $M$, every strongly prime
        element of $(\Q(\Sig M),\sqsubseteq)$ is principal.
  \item If $X$ is a well-filtered space, then $\up x$ is strongly prime in $(\Q(X),\sqsubseteq)$
        for every $x\in X$.
\end{enumerate}
\end{lemma}

\begin{proof}
(i) Let $T\in\Q(\Sig M)$ be strongly prime; we write $T$ as the meet of a Scott-compact family. The
map $\xi\colon M\to\Q(\Sig M)$, $\xi(m)=\up m$,
is monotone, because $m\le m'$ implies $\up m\supseteq\up m'$. If $E\subseteq M$ is directed and has
a supremum $s$, then $\up s=\bigcap_{e\in E}\up e$ is the supremum of $\xi(E)$ in
$(\Q(\Sig M),\sqsubseteq)$: it is an upper bound, and every upper bound $Q$ satisfies $Q\subseteq\up
e$ for all $e\in E$, hence $Q\subseteq\up s$. So $\xi$ preserves the existing directed suprema, and
therefore it is continuous for the Scott topologies of $M$ and of $(\Q(\Sig M),\sqsubseteq)$
(Section~\ref{sec:prelim}). Hence
$\xi(T)$ is compact in the latter topology, and so is its saturation
$\mathcal K=\{Q\in\Q(\Sig M): Q\subseteq\up m\text{ for some }m\in T\}$.
Every member of $\mathcal K$ is included in $T$, because $T$ is an upper set; and a set
$Q\in\Q(\Sig M)$ that includes all members of $\mathcal K$ includes $\bigcup_{m\in T}\up m=T$. So
$T$ is the meet of $\mathcal K$. As $T$ is strongly prime, $T\in\mathcal K$: $T\subseteq\up m$ for
some $m\in T$, so $T=\up m$.

(ii) Let $\mathcal K$ be as in the definition, with $\bigwedge\mathcal K\sqsubseteq\up x$. We show
that $\bigcup\mathcal K$ is compact; then it is the meet of $\mathcal K$, and $x$ lies in a member
of $\mathcal K$. As $X$ is well-filtered, the supremum of a directed subset $\mathcal H$ of
$(\Q(X),\sqsubseteq)$ is the intersection $\bigcap\mathcal H$, which belongs to $\Q(X)$ by
Section~\ref{sec:prelim}. For every open set $U$ of $X$ the set $\Box U=\{Q\in\Q(X):Q\subseteq
U\}$ is Scott-open in $(\Q(X),\sqsubseteq)$, because it is an upper set and $\bigcap\mathcal
H\subseteq U$ implies that some member of $\mathcal H$ is included in $U$. Let $\mathcal O$ be a family of
open subsets of $X$, closed under finite unions, that covers $\bigcup\mathcal K$. Each member of
$\mathcal K$ is compact, hence included in a member of $\mathcal O$, as $\mathcal O$ is closed under
finite unions. So the sets $\Box O$, $O\in\mathcal O$, cover $\mathcal K$, and finitely many of
them cover $\mathcal K$; hence $\bigcup\mathcal K$ is included in one member of $\mathcal O$. By
the criterion for compactness in Section~\ref{sec:prelim}, $\bigcup\mathcal K$ is compact; it is
saturated, so it is the meet $\bigwedge\mathcal K$. As $\up x\subseteq\bigcup\mathcal K$,
some $Q\in\mathcal K$ contains $x$. Then $\up x\subseteq Q$, that is, $Q\sqsubseteq\up x$, and
$\up x\in\mathcal K$ because $\mathcal K$ is an upper set.
\end{proof}

\begin{proof}[Proof of Proposition~\ref{prop:HL}]
Prime and strongly prime elements are defined by the order alone, so order isomorphisms preserve
them. Let $L$ be a poset such that $\Sig L$ is locally compact and well-filtered, let $M$ be a
poset, and let $F$ be an order isomorphism of $(\Q(\Sig L),\sqsubseteq)$ onto
$(\Q(\Sig M),\sqsubseteq)$.

\emph{Step 1: $F(\up x)$ is principal for every $x\in L$.} The set $\up x$ is strongly prime by
Lemma~\ref{lem:sprime}(ii); so is $F(\up x)$, which is principal by Lemma~\ref{lem:sprime}(i).

\emph{Step 2: $F^{-1}(\up m)$ is principal for every $m\in M$.} The set $\up m$ is prime, so
$F^{-1}(\up m)$ is prime, and it is principal by Lemma~\ref{lem:prime}, as $\Sig L$ is locally
compact.

\emph{Step 3: $L$ and $M$ are isomorphic.} By Steps~1 and~2, $F$ maps $\{\up x:x\in L\}$ onto
$\{\up m:m\in M\}$, and $\up x\sqsubseteq\up y$ if and only if $x\le y$, in $L$ and in $M$.
\end{proof}

\begin{corollary}\label{cor:HL}
The dcpo $P$ of Theorem~\ref{thm:domain} is $\Q_\sigma$-unique, although $\Sig P$ is not co-sober.
So co-sobriety of the Scott space is not necessary for a poset to be $\Q_\sigma$-unique.
\end{corollary}

\begin{proof}
By Theorem~\ref{thm:domain}, $\Sig P$ is not co-sober, and it is locally compact and sober, hence
well-filtered (Section~\ref{sec:prelim}); so Proposition~\ref{prop:HL} applies.
\end{proof}

\subsection{Sobrification}
A \emph{sobrification} of a $T_0$ space $Y$ is a sober space $Y^s$ together with a continuous map
$j\colon Y\to Y^s$ such that $O\mapsto j^{-1}(O)$ is a bijection from the open sets of $Y^s$ onto
those of $Y$; it is unique up to homeomorphism \cite[Chapter~8]{GL13}.

\begin{corollary}\label{cor:sobrification}
Every poset with its Alexandroff topology is co-sober \textup{(\cite[Lemma~2.1(2)]{JX26sp})}. The space $\Sig P$ of
Theorem~\ref{thm:domain} is a sobrification of the poset $(\B,\supseteq)$ with its Alexandroff
topology; this poset is isomorphic to the poset of compact elements of $P$. So sobrification does
not preserve co-sobriety.
\end{corollary}

\begin{proof}
We recall the proof of the first statement. Let $Y$ be a poset with its Alexandroff topology; its
specialization order is the order of $Y$. A
compact saturated set $H$ is covered by the open sets $\up q$, $q\in H$, hence by finitely many of
them; so $H$ is a finite union of principal sets, and if $H$ is k-irreducible, it is one of them.

Now let $Y=(\B,\supseteq)$ and $j\colon Y\to\Sig P$, $j(U)=\dn U$. As $\dn U\subseteq\dn V$ if and
only if $V\subseteq U$, $j$ is an isomorphism of $Y$ onto its image, the poset of compact elements of
$P$ (Lemma~\ref{lem:idl}(i)). For $U\in\B$ we have
$j^{-1}([U])=\{V\in\B:V\subseteq U\}=\up_YU$. So $j$ is continuous, and every upper set $\mathcal U$
of $Y$ is $j^{-1}\bigl(\bigcup_{U\in\mathcal U}[U]\bigr)$. If $O$ and $O'$ are Scott-open sets with
$j^{-1}(O)=j^{-1}(O')$, then $O$ and $O'$ contain the same compact elements of $P$, so $O=O'$ by
Lemma~\ref{lem:idl}(i) and~(ii). As $\Sig P$ is sober, it is a sobrification of $Y$.
\end{proof}

\begin{remark}\label{rem:domain}
\textup{(a)} \emph{$K$ is not prime}, in accordance with Lemma~\ref{lem:prime}. Let $p\ne q$ be
points of $\N$. The sets $[\N\setminus\{p\}]$ and $[\N\setminus\{q\}]$ (Lemma~\ref{lem:base}(i)) are
compact and saturated. Their union includes $K$, because every $m\in\N$ differs from $p$ or from $q$.
Neither of them
includes $K$, because $I_p\notin[\N\setminus\{p\}]$ and $I_q\notin[\N\setminus\{q\}]$. So
k-irreducible sets need not be prime, even in the Scott space of an algebraic dcpo.

\textup{(b)} \emph{Coherence fails for two compact open sets, $\B$ is not closed under finite
intersections, and $P$ is not a lattice.} Let $a\in\A$ and $x=\ell(a)$; then
$O_a,\N\setminus\{x\}\in\B$ (Lemma~\ref{lem:base}(i)). The set $G=[O_a]\cap[\N\setminus\{x\}]$ is
the union of the sets $[W]$ with $W\in\B$ and $W\subseteq O_a\setminus\{x\}$. Each such $W$ meets
$a$ in a finite set by $(\mathrm R_a)$, and $I_d\in[W]$ only if $d\in W$. But $I_d\in G$ for every
$d\in a\setminus(Z_a\cup\{x\})$, and $a\cap Z_a\subseteq(a\cap c(a))\cup\{\ell(c(a))\}$ is finite.
So no finite number of the sets $[W]$ covers $G$; in particular $O_a\setminus\{x\}\notin\B$. As
$G$ is the set of common upper bounds of $\dn O_a$ and $\dn(\N\setminus\{x\})$ in $P$, these have
no supremum $z$; otherwise $G$ would be the compact set $\up z$.

\textup{(c)} \emph{In a sober space, if $Q$ is compact and saturated and $\min Q$ is irreducible and
not a singleton, then $\min Q$ is not closed in $Q$. This applies to $K$, as $\min K=e(\N)$ is
homeomorphic to the anti-Hausdorff, hence irreducible, space $D$ (Section~\ref{sec:prelim}).} By
sobriety the closure of $\min Q$ is $\dn z$ for some $z$, an upper bound of the antichain $\min Q$,
which is not a singleton; so $z\in Q\setminus\min Q$.

\textup{(d)} \emph{The plan of \S\ref{ssec:idea} fails for the base $\B'$ of all nonempty open
sets.} By the proof of Lemma~\ref{lem:idl}, the Scott space of the ideal completion
$P'$ of $(\B',\supseteq)$ is locally compact and sober, hence well-filtered. Its basic compact open
sets $[U]$, $U\in\B'$, are closed under finite intersections, as $\B'$ is ($D$ is anti-Hausdorff). As
they form a base, each compact saturated set is, by compactness, the intersection of the filtered
family of finite unions of basic sets that include it. So the intersection of two compact
saturated sets is a filtered intersection of compact open sets, hence compact
(Section~\ref{sec:prelim}). Thus $\Sig P'$ is coherent, hence co-sober \cite{ZHW23}, and
the analogue of $K$ is not k-irreducible.

\textup{(e)} \emph{$\Sig P$ is not first-countable, and the method cannot give a first-countable
example.} Suppose that $X$ is first-countable, $Q\in\Q(X)$, and $Q'\cap\min Q$ is closed in
$\min Q$ for every compact set $Q'\subseteq Q$, as in Lemma~\ref{lem:dtrace}. Then every compact
subset of $\min Q$ is closed in $\min Q$; that is, $\min Q$ is a first-countable KC space, hence
Hausdorff, as is classical \cite{Wil67}. But the method needs $\min Q$ to be anti-Hausdorff with at
least two points. In particular $\Sig P$ is not first-countable: by Lemma~\ref{lem:dtrace} the
hypothesis holds for $K$, and $\min K=e(\N)$ is a copy of the anti-Hausdorff space $D$.
\end{remark}

\begin{question}\label{q:firstcountable}
Is every first-countable, locally compact, sober space co-sober? In particular, is the Scott space
of a continuous dcpo co-sober if it is first-countable?
\end{question}

\subsection*{Use of artificial intelligence}

During the preparation of this work the authors used Claude (Anthropic), an AI assistant, in order
to search the literature, to develop and check the arguments, and to draft the text. After using
this tool, the authors reviewed and edited the content as needed and take full responsibility for
the content of the article.


\begin{thebibliography}{99}
\bibitem{ELS04} M.~Escard\'o, J.~Lawson, A.~Simpson, \emph{Comparing Cartesian closed categories of
  (core) compactly generated spaces}, Topology Appl.\ \textbf{143} (2004), 105--145.
\bibitem{Giam26} K.~K.~Giam, \emph{Retracts of well-filtered co-sober spaces are co-sober}, Topology
  Appl., in press (2026).
\bibitem{GHKLMS03} G.~Gierz, K.~H.~Hofmann, K.~Keimel, J.~D.~Lawson, M.~Mislove, D.~S.~Scott,
  \emph{Continuous lattices and domains}, Encyclopedia of Mathematics and its Applications, vol.~93,
  Cambridge University Press, Cambridge, 2003.
\bibitem{GL13} J.~Goubault-Larrecq, \emph{Non-Hausdorff topology and domain theory}, New
  Mathematical Monographs, vol.~22, Cambridge University Press, Cambridge, 2013.
\bibitem{HL25} H.~Hou, Q.~Li, \emph{Posets uniquely determined by its Scott compact saturated
  subsets}, preprint, arXiv:2503.02602, 2025.
\bibitem{Hru14} M.~Hru\v{s}\'ak, \emph{Almost disjoint families and topology}, Recent progress in
  general topology III (K.~P.~Hart, J.~van Mill, P.~Simon, eds.), Atlantis Press, Paris, 2014,
  pp.~601--638.
\bibitem{JX26sp} W.~Ji, X.~Xu, \emph{A co-sober Alexandrov space whose Scott power space is not
  co-sober}, preprint, arXiv:2611.01516, 2026.
\bibitem{JX26cba} W.~Ji, X.~Xu, \emph{Scott spaces of complete Boolean algebras need not be
  co-sober}, preprint, arXiv:2608.08484, 2026.
\bibitem{Kov05} M.~M.~Kov\'ar, \emph{The Hofmann--Mislove theorem for general topological
  structures}, Spatial representation: discrete vs.\ continuous computational models
  (R.~Kopperman, M.~B.~Smyth, D.~Spreen, J.~Webster, eds.), Dagstuhl Seminar Proceedings, vol.~04351,
  Schloss Dagstuhl -- Leibniz-Zentrum f\"ur Informatik, 2005, pp.~1--9.
\bibitem{SWXZ20} C.~Shen, G.~Wu, X.~Xi, D.~Zhao, \emph{Sober Scott spaces are not always co-sober},
  Topology Appl.\ \textbf{282} (2020), 107316.
\bibitem{Sim80} P.~Simon, \emph{A compact Fr\'echet space whose square is not Fr\'echet}, Comment.\
  Math.\ Univ.\ Carolin.\ \textbf{21} (1980), 749--753.
\bibitem{vD93} E.~K.~van Douwen, \emph{An anti-Hausdorff Fr\'echet space in which convergent
  sequences have unique limits}, Topology Appl.\ \textbf{51} (1993), 147--158.
\bibitem{WX18} X.~Wen, X.~Xu, \emph{Sober is not always co-sober}, Topology Appl.\ \textbf{250}
  (2018), 48--52.
\bibitem{Wil67} A.~Wilansky, \emph{Between $T_1$ and $T_2$}, Amer.\ Math.\ Monthly \textbf{74}
  (1967), 261--266.
\bibitem{XXJQ26k} X.~Xu, F.~Xu, H.~Ji, D.~Qiu, \emph{Kernel retracts of co-sober spaces need not be
  co-sober}, Topology Appl., in press; preprint arXiv:2611.01531 (2026), under the title
  \emph{Kernel retracts of co-sober spaces}.
\bibitem{XXJQ26s} X.~Xu, F.~Xu, H.~Ji, D.~Qiu, \emph{Smyth power space of a co-sober space is not
  always co-sober}, preprint, arXiv:2608.02047, 2026.
\bibitem{ZHW23} B.~Zhao, Z.~He, K.~Wang, \emph{On co-sober spaces}, Houston J.~Math.\ \textbf{49}
  (2023), no.~4, 971--988.
\end{thebibliography}
\end{document}